\def\titlep{Tensor products of type III factor 
representations of Cuntz-Krieger algebras}
\documentclass[11pt]{article}
\usepackage{graphicx,ifthen}
\usepackage{amssymb}
\usepackage{amsmath}

\font\germ=eufm10 at12pt

\def\goth#1{\hbox{\germ#1}}

\newcommand{\qed}{\hbox{\rule[-2pt]{3pt}{6pt}}}
\newcommand{\qedh}{\hfill\qed \\}

\newcommand{\vv}{\vspace{.3in}}

\newtheorem{Thm}{Theorem}[section]

\newtheorem{rem}[Thm]{Remark}

\newtheorem{ex}[Thm]{Example}
\newtheorem{defi}[Thm]{Definition}
\newtheorem{lem}[Thm]{Lemma}

\newtheorem{prop}[Thm]{Proposition}

\newtheorem{cor}[Thm]{Corollary}

\newcommand{\ww}{\vv\noindent}

\def\cal#1{\mathcal #1}
\def\con{{\cal O}_{n}}

\def\pr{{\it Proof.}\quad}

\def\nset#1{\{1,\ldots,n\}^{#1}}

\def\coa{{\cal O}_{A}}
\def\co#1{{\cal O}_{#1}}
\def\disp#1{{\displaystyle #1}}
\def\ck{Cuntz-Krieger}

\def\brl{branching law}
\def\bfsnl{{\rm BFS}_{N}(\Lambda)}

\addtocounter{footnote}{1}
\def\sftt#1{
\setcounter{equation}{0}
\addtocounter{footnote}{1}
\section{#1}
}

\def\ssft#1{\subsection{#1}}
\def\ssfr#1{\subsection*{#1}}

\begin{document}
%
%
\def\autherp{Katsunori Kawamura}
\def\emailp{e-mail: kawamura@kurims.kyoto-u.ac.jp.}
\def\addressp{{\small {\it College of Science 
and Engineering Ritsumeikan University,}}\\
{\small {\it 1-1-1 Noji Higashi, Kusatsu, Shiga 525-8577, Japan}}
}

\def\infw{\Lambda^{\frac{\infty}{2}}V}
\def\zhalfs{{\bf Z}+\frac{1}{2}}
\def\ems{\emptyset}
\def\pmvac{|{\rm vac}\!\!>\!\! _{\pm}}
\def\vac{|{\rm vac}\rangle _{+}}
\def\dvac{|{\rm vac}\rangle _{-}}
\def\ovac{|0\rangle}
\def\tovac{|\tilde{0}\rangle}
\def\expt#1{\langle #1\rangle}
\def\zph{{\bf Z}_{+/2}}
\def\zmh{{\bf Z}_{-/2}}
\def\brl{branching law}
\def\bfsnl{{\rm BFS}_{N}(\Lambda)}
\def\scm#1{S({\bf C}^{N})^{\otimes #1}}
\def\mqb{\{(M_{i},q_{i},B_{i})\}_{i=1}^{N}}
\def\zhalf{\mbox{${\bf Z}+\frac{1}{2}$}}
\def\zmha{\mbox{${\bf Z}_{\leq 0}-\frac{1}{2}$}}
\newcommand{\mline}{\noindent
\thicklines
\setlength{\unitlength}{.1mm}
\begin{picture}(1000,5)
\put(0,0){\line(1,0){1250}}
\end{picture}
\par
 }
\def\ptimes{\otimes_{\varphi}}
\def\delp{\Delta_{\varphi}}
\def\delps{\Delta_{\varphi^{*}}}
\def\gamp{\Gamma_{\varphi}}
\def\gamps{\Gamma_{\varphi^{*}}}
\def\sem{{\sf M}}
\def\hdelp{\hat{\Delta}_{\varphi}}
\def\tilco#1{\tilde{\co{#1}}}
\def\ndm#1{{\bf M}_{#1}(\{0,1\})}
\def\cdm#1{{\cal M}_{#1}(\{0,1\})}
\def\tndm#1{\tilde{{\bf M}}_{#1}(\{0,1\})}
\def\sck{{\sf CK}_{*}}
\def\hdel{\hat{\Delta}}
\def\ba{\mbox{\boldmath$a$}}
\def\bb{\mbox{\boldmath$b$}}
\def\bc{\mbox{\boldmath$c$}}
\def\be{\mbox{\boldmath$e$}}
\def\bp{\mbox{\boldmath$p$}}
\def\bq{\mbox{\boldmath$q$}}
\def\bu{\mbox{\boldmath$u$}}
\def\bv{\mbox{\boldmath$v$}}
\def\bw{\mbox{\boldmath$w$}}
\def\bx{\mbox{\boldmath$x$}}
\def\by{\mbox{\boldmath$y$}}
\def\bz{\mbox{\boldmath$z$}}
\def\bomega{\mbox{\boldmath$\omega$}}
\def\N{{\bf N}}
\def\lxm{L_{2}(X,\mu)}
%
%
%
\setcounter{section}{0}
\setcounter{footnote}{0}
\setcounter{page}{1}
\pagestyle{plain}

%
%
\title{\titlep}
\author{\autherp\thanks{\emailp}
\\
\addressp}
\date{}
\maketitle

%
%
\begin{abstract}
We introduced a non-symmetric tensor product 
of any two states or any two representations 
of Cuntz-Krieger algebras associated with
a certain non-cocommutative comultiplication in previous our work.
In this paper,
we show that a certain set of KMS states is closed with respect 
to the tensor product.
From this, we obtain formulae of tensor product 
of type {\rm III} factor representations of Cuntz-Krieger algebras
which is different from results of the tensor product 
of factors of type {\rm III}.
\end{abstract}

\noindent
{\bf Mathematics Subject Classifications (2000).} 
46K10, 46L30, 46L35.
\\
{\bf Keywords.} tensor product of representations, Cuntz-Krieger algebra,
KMS state, type {\rm III} factor representation.
%
%
\sftt{Introduction}
\label{section:first}
We have studied states and representations of operator algebras.
KMS states over Cuntz-Krieger algebras with respect to certain one-parameter 
automorphism groups are known \cite{EFW,EL2}.
GNS representations by such KMS states 
induce type III factor representations of Cuntz-Krieger algebras \cite{Okayasu}.
On the other hand, we introduced a non-symmetric tensor product 
of any two states or any two representations of 
Cuntz-Krieger algebras associated with
a certain non-cocommutative comultiplication \cite{TS05}.
From these,
we investigate the tensor product of two KMS states or their GNS representations
of Cuntz-Krieger algebras in this paper.
In consequence, we obtain formulae of tensor product 
of type {\rm III} factor representations of Cuntz-Krieger algebras.
Remark that this is not a study of the tensor product of 
``factors" of type {\rm III} but that of 
a certain non-symmetric tensor product of 
 type {\rm III} ``factor representations."

%
%
\ssft{Motivation}
\label{subsection:firstone} 
From studies of branching laws of a certain class of representations 
of Cuntz algebras
\cite{PE01,PE02,PE03},
we found a non-symmetric tensor product of representations 
of Cuntz algebras \cite{TS01}. 
Subsequently,
this tensor product was explained by the C$^{*}$-bialgebra which is defined by
the direct sum of all Cuntz algebras \cite{TS02}.
Furthermore such a construction of C$^{*}$-bialgebra was generalized to 
Cuntz-Krieger algebras \cite{TS05}.
The essential tool of the construction is a certain set of embeddings
among Cuntz-Krieger algebras.
We explain this as follows.

A matrix $A$ is {\it nondegenerate} if any column and any row are not zero.
For $2\leq n<\infty$,
let $\ndm{n}$ denote the set of all
nondegenerate $n\times n$ matrices with entries $0$ or $1$
and define ${\bf M}\equiv \cup_{n}\,\ndm{n}$.
For $A\in {\bf M}$,
let $\coa$ denote the Cuntz-Krieger algebra by $A$ \cite{CK}.
For $A,B\in {\bf M}$,
let $A\boxtimes B$ denote the Kronecker product of $A$ and $B$ \cite{Dief}.
Then ${\bf M}$ is closed with respect to $\boxtimes$.
We can construct a set $\{\varphi_{A,B}:A,B\in {\bf M}\}$
such that $\varphi_{A,B}$  is 
a unital $*$-embedding of $\co{A\boxtimes B}$ into $\coa\otimes\co{B}$ and 
%
%
\begin{equation}
\label{eqn:weak}
(\varphi_{A,B}\otimes id_{C})\circ \varphi_{A\boxtimes B,C}=
(id_{A}\otimes \varphi_{B,C})\circ \varphi_{A, B\boxtimes C}
\quad(A,B,C\in {\bf M})
\end{equation}
where $id_{X}$ denotes the identity map of $\co{X}$ for $X=A,C$ 
and $\coa\otimes \co{B}$ means the minimal tensor product 
of $\coa$ and $\co{B}$.
We will give the explicit definition of these embeddings 
in $\S$ \ref{subsection:firsttwo}.

From the set $\varphi\equiv \{\varphi_{A,B}:A,B\in {\bf M}\}$,
we can define associative tensor products
of states or representations as follows.
Let ${\cal S}_{A}$ denote the set of all states over $\coa$.
For $\rho_{1}\in {\cal S}_{A}$
and $\rho_{2}\in {\cal S}_{B}$,
define $\rho_{1}\ptimes \rho_{2}\in  {\cal S}_{A\boxtimes B}$ by
%
%
\begin{equation}
\label{eqn:tensorstate}
\rho_{1}\ptimes \rho_{2}
\equiv (\rho_{1}\otimes \rho_{2})\circ \varphi_{A,B}.
\end{equation}
From (\ref{eqn:weak}), we see that
\[(\rho_{1}\ptimes \rho_{2})\ptimes \rho_{3}=
\rho_{1}\ptimes (\rho_{2}\ptimes \rho_{3})\quad
(\rho_{1},\rho_{2},\rho_{3}\in {\cal S}_{*})\]
where
${\cal S}_{*}\equiv \bigcup_{A\in {\bf M}}{\cal S}_{A}$.
In consequence,
${\cal S}_{*}$ is a semigroup with respect to the operation $\ptimes$.
Let ${\rm Rep}\coa$ denote the class of all unital $*$-representations 
of $\coa$.
As is the case with states,
we can define the associative operation
$\ptimes$ on the class ${\sf R}_{*}\equiv \bigcup_{A\in {\bf M}}{\rm Rep}\coa$:
\[\ptimes: {\sf R}_{*}\times {\sf R}_{*}\to {\sf R}_{*}.\]
We show properties of $\ptimes$ as follows:
For two representations
$\pi_{1},\pi_{2}$ of a C$^{*}$-algebra ${\cal A}$,
if $\pi_{1}$ and $\pi_{2}$ are unitarily equivalent
({\it resp.} quasi-equivalent), then 
we write $\pi_{1}\simeq\pi_{2}$ ({\it resp.} $\pi_{1}\approx \pi_{2}$).
\begin{enumerate}
\item
For $\pi_{i},\pi_{i}^{'}\in {\sf R}_{*}$,
if $\pi_{i}\simeq\pi_{i}^{'}$ for each $i=1,2$,
then $\pi_{1}\ptimes \pi_{2}\simeq
\pi_{1}^{'}\ptimes \pi_{2}^{'}$.
From this,
$\ptimes$ is well-defined on the set 
%
%
\begin{equation}
\label{eqn:rtwo}
{\cal R}_{*}\equiv 
\bigcup_{A\in\cdm{*}}({\rm Rep}\coa/\!\!\simeq)
\end{equation}
of all unitary equivalence classes of representations of $\coa$'s.
\item
The operation 
$\ptimes$ satisfies the distribution law with respect to the direct sum. 
\item
The operation $\ptimes$ is non-symmetric in a strong sense,
that is, there exist $A\in {\bf M}$ and $\pi_{1}, \pi_{2}\in {\rm Rep}\coa$
such that $\pi_{1}\ptimes \pi_{2}\not\simeq\pi_{2}\ptimes \pi_{1}$.
\item
Even if both $\pi_{1}$ and $\pi_{2}$ are irreducible,
$\pi_{1}\ptimes \pi_{2}$ is not always irreducible ($\S$ 4.1 of \cite{TS01}).
\end{enumerate}
In addition, we can show the following.
%
%
\begin{lem}
\label{lem:factor}
\begin{enumerate}
\item
There exist $A\in {\bf M}$ and 
two factor representations $\pi_{1}$ and $\pi_{2}$ of $\coa$
such that $\pi_{1}\ptimes \pi_{2}$ is not a factor representation
of $\co{A\boxtimes A}$.
\item
There exist $A\in {\bf M}$ and $\pi_{1}, \pi_{2}\in {\rm Rep}\coa$
such that $\pi_{1}\ptimes \pi_{2}\not\approx\pi_{2}\ptimes \pi_{1}$.
\item
For $\pi_{i},\pi_{i}^{'}\in {\sf R}_{*}$,
if $\pi_{i}\approx \pi_{i}^{'}$ for each $i=1,2$,
then $\pi_{1}\ptimes \pi_{2}\approx\pi_{1}^{'}\ptimes \pi_{2}^{'}$.
\item
There exist $\rho_{1},\rho_{2}\in {\cal S}_{*}$ 
such that 
$\pi_{\rho_{1}}\ptimes \pi_{\rho_{2}}\not\approx\pi_{\rho_{1}\ptimes \rho_{2}}$ 
where $\pi_{\rho}$ denotes the GNS representation by a state $\rho$.
\end{enumerate}
\end{lem}

\noindent
From Lemma \ref{lem:factor}(i),
even if both $\pi_{1}$ and $\pi_{2}$ are factor representations,
$\pi_{1}\ptimes \pi_{2}$ is not always a factor representation.
From Lemma \ref{lem:factor}(ii),
$\ptimes$ is also non-symmetric with respect to quasi-equivalence classes.
Lemma \ref{lem:factor} will be proved in $\S$ \ref{subsection:thirdfour}.

Our interest is to compute tensor products of 
concrete states or representations of
Cuntz-Krieger algebras with respect to the above operation $\ptimes$.
For the case of permutative representations of Cuntz algebras,
we gave formulae of decomposition of tensor products completely \cite{TS01}. 
In this case, non-type I representation never appear.
In this paper,
we intend to consider tensor products of non-type I representations.

%
%
\ssft{A set of embeddings of Cuntz-Krieger algebras}
\label{subsection:firsttwo}
In this subsection, we review a set of embeddings of 
Cuntz-Krieger algebras \cite{TS05}.
We state that a matrix $A\in M_{n}({\bf C})$ is {\it irreducible} if 
for any $i,j\in\{1,\ldots,n\}$,
there exists $k\in {\bf N}\equiv \{1,2,3,\ldots\}$ such that
$(A^{k})_{i,j}\ne 0$ where $A^{k}=A\cdots A$ ($k$ times).
For $2\leq n<\infty$,
let $\cdm{n}$ denote the set of all 
irreducible nondegenerate $n\times n$ matrices with entries $0$ or $1$,
which is not a permutation matrix.
Define
%
%
\begin{equation}
\label{eqn:nondegenerate}
\cdm{*}\equiv \cup \{\cdm{n}:n\geq 2\}.
\end{equation}

For $A=(A_{ij})\in M_{n}({\bf C})$ and $B=(B_{ij})\in M_{m}({\bf C})$,
define the {\it Kronecker product} $A\boxtimes B\in M_{nm}({\bf C})$ 
of $A$ and $B$ by
%
%
\begin{equation}
\label{eqn:phitwo}
(A\boxtimes B)_{m(i-1)+j,m(i^{'}-1)+j^{'}}\equiv A_{ii^{'}}B_{jj^{'}}
\end{equation}
for $i,i^{'}\in \nset{}$ and $j,j^{'}\in\{1,\ldots,m\}$ \cite{Dief}.
If $A,B\in\cdm{*}$, then $A\boxtimes B\in \cdm{*}$.

For $A\in\cdm{n}$,
let $s_{1}^{(A)},\ldots,s_{n}^{(A)}$ 
denote the canonical generators of $\coa$.
For $A\in\cdm{n}$ and $B\in\cdm{m}$,
define the map $\varphi_{A,B}$ from $\co{A\boxtimes B}$ 
to the minimal tensor product $\co{A}\otimes \co{B}$ by
%
%
\begin{equation}
\label{eqn:varphi}
\varphi_{A,B}(s_{m(i-1)+j}^{(A\boxtimes B)})
\equiv s_{i}^{(A)}\otimes s_{j}^{(B)}\quad 
(i\in\nset{},\,j\in\{1,\ldots,m\}).
\end{equation}
Then we can verify that $\varphi_{A,B}$ 
is uniquely extended to a unital $*$-embedding of 
$\co{A\boxtimes B}$ into $\co{A}\otimes \co{B}$.
Then $\{\varphi_{A,B}:A,B\in\cdm{*}\}$ satisfies (\ref{eqn:weak}).

%
%
\ssft{GNS representations of Cuntz-Krieger algebras 
by certain states}
\label{subsection:firstthree}
In this subsection, we introduce an 
index set of states and representations of Cuntz-Krieger algebras
by \cite{EL,Okayasu}
such that it is suitable to compute the tensor product by $\ptimes$.
We rewrite original statements as follows.
Define $I_{0}\equiv \{x\in {\bf R}:0<x<1\}$.
For $\ba=(a_{1},\ldots,a_{n})\in I_{0}^{n}$,
define $\hat{\ba}\equiv {\rm diag}(a_{1},\ldots,a_{n})\in M_{n}({\bf R})$.
For $A\in \cdm{n}$,
let $\hat{\ba}A$ denote the product of matrices $\hat{\ba}$ and $A$.
Define
%
%
\begin{equation}
\label{eqn:pfe}
\Lambda(A)\equiv \{\ba\in I_{0}^{n}:PFE(\hat{\ba}A)=1\}
\end{equation}
where $PFE(X)$ denotes the Perron-Frobenius eigenvalue 
of a irreducible non-negative matrix $X$ \cite{BP,Seneta}. 
Since 
%
%
\begin{equation}
\label{eqn:eigentwo}
\be(A)\equiv (1/c_{A},\ldots,1/c_{A})
\end{equation}
belongs to $\Lambda(A)$ for $c_{A}\equiv PFE(A)$,
$\Lambda(A)\ne \emptyset$ for each $A\in\cdm{*}$.

Define ${\bf R}_{+}\equiv \{x\in {\bf R}:x>0\}$.
For $A\in \cdm{n}$ and $\ba=(a_{1},\ldots,a_{n})\in \Lambda(A)$,
let $\bx=(x_{1},\ldots,x_{n})\in {\bf R}^{n}_{+}$ 
denote the Perron-Frobenius eigenvector 
of $\hat{\ba}A$ such that $x_{1}+\cdots+x_{n}=1$ and
let $s_{1},\ldots,s_{n}$ denote canonical generators of $\coa$.
Define the state $\rho_{\ba}$ over $\coa$ by
%
%
\begin{equation}
\label{eqn:kmsthree}
\rho_{\ba}(s_{J}s_{K}^{*})
=\delta_{JK}a_{j_{1}}\ldots a_{j_{m-1}}x_{j_{m}}
\end{equation}
when $s_{J}s_{K}^{*}\ne 0$
for $J=(j_{1},\ldots,j_{m})\in\{1,\ldots,n\}^{m}$
and $K\in \cup_{l\geq 1}\{1,\ldots,n\}^{l}$
where $s_{J}=s_{j_{1}}\cdots s_{j_{m}}$.
Then  the following holds.
%
%
\begin{Thm}
\label{Thm:typetwo}
For $\ba\in\Lambda(A)$,
let $\varpi_{\ba}$ denote the GNS representation of $\coa$ by $\rho_{\ba}$.
Then the von Neumann algebra $M_{\ba}\equiv \varpi_{\ba}(\coa)^{''}$ is an 
approximately finite dimensional (=AFD) factor of type ${\rm III}$.
Furthermore,
Connes' classification of the type of $M_{\ba}$ \cite{Connes}
is given as follows:
\begin{enumerate}
\item
If there exist $p_{1},\ldots,p_{n}\in {\bf N}$ and $0<\lambda<1$
such that  the greatest common divisor of the set $\{p_{1},\ldots,p_{n}\}$
is $1$ and $\ba=(\lambda^{p_{1}},\ldots,\lambda^{p_{n}})$,
then $p_{1},\ldots,p_{n},\lambda$ are uniquely determined by $\ba$,
and $M_{\ba}$ is of type ${\rm III}_{\lambda}$.
\item
If assumptions in (i) do not hold,
then 
$M_{\ba}$ is of type ${\rm III}_{1}$.
\end{enumerate}
\end{Thm}
In addition to Theorem \ref{Thm:typetwo},
we define the positive real number $\lambda(\ba)$ by
%
%
\begin{equation}
\label{eqn:lambdatwo}
\lambda(\ba)\equiv \left\{
\begin{array}{ll}
\lambda \quad&( \ba \mbox{ is as in the case of (i)}),\\
\\
1 \quad&( \ba \mbox{ is as in the case of (ii)}).\\
\end{array}
\right.
\end{equation}
Then we can simply state that $M_{\ba}$ is of type ${\rm III}_{\lambda(\ba)}$
for each $\ba\in \Lambda(A)$.
Especially,
for $\be(A)$ in (\ref{eqn:eigentwo}), $\lambda(\be(A))=1/c_{A}$.
The condition in Theorem \ref{Thm:typetwo}(ii) is easily rewritten as follows.
%
%
\begin{lem}
\label{lem:log}
For $\ba=(a_{1},\ldots,a_{n})\in\Lambda(A)$,
$\lambda(\ba)=1$
if and only if 
there exist $i,j$ such that
$\log a_{i}/\log a_{j}\not\in {\bf Q}$.
\end{lem} 
Theorem \ref{Thm:typetwo} is nothing but a reformulation 
of Theorem 4.2 in \cite{Okayasu}
without use of terminology of KMS states.
Here we abbreviate one-parameter automorphism groups and inverse temperatures
associated with above KMS states for simplicity of description.
These will be  explained in $\S$ \ref{section:third}.
%
%
\begin{rem}
\label{rem:lambda}
{\rm
\begin{enumerate}
\item
For $A=(A_{ij})$,
assume $A_{ij}=1$ for each $i,j$.
By definition,
$\lambda(a_{1},\ldots,a_{n})=\lambda(a_{p(1)},\ldots,a_{p(n)})$
for any $(a_{1},\ldots,a_{n})\in\Lambda(A)$ 
and any permutation $p\in{\goth S}_{n}$.
Therefore $\ba\mapsto \lambda(\ba)$ is not injective in general.
\item
For $\ba,\bb\in \Lambda(A)$,
if $\lambda(\ba)\ne \lambda(\bb)$,
then $\varpi_{\ba}$ and $\varpi_{\bb}$ 
are not quasi-equivalent.
However,
we do not know whether
$\varpi_{\ba}$ and $\varpi_{\bb}$ 
are quasi-equivalent or not
even if $\lambda(\ba)=\lambda(\bb)$.
\end{enumerate}
}
\end{rem}

%
%
\ssft{Main theorems}
\label{subsection:firstfour}
In this subsection, we show our main theorems.
For vectors $\bv=(v_{1},\ldots,v_{n})\in {\bf C}^{n}$ 
and $\bw=(w_{1},\ldots,w_{m})\in {\bf C}^{m}$,
define $\bv\boxtimes \bw\in {\bf C}^{nm}$ by
%
%
\begin{equation}
\label{eqn:vector}
(\bv\boxtimes \bw)_{m(i-1)+j}=v_{i}w_{j}\quad
(i=1,\ldots,n,\,j=1,\ldots,m).
\end{equation}
We see that $(\bv\boxtimes \bw)\boxtimes \bx=\bv\boxtimes (\bw\boxtimes \bx)$.
For $\cdm{*}$ in (\ref{eqn:nondegenerate}) and $\Lambda(A)$ in (\ref{eqn:pfe}),
define
%
%
\begin{equation}
\label{eqn:lambdazero}
\Lambda(*)\equiv \bigcup_{A\in\cdm{*}}\Lambda(A).
\end{equation}
Then the following holds for $\ptimes$ in (\ref{eqn:tensorstate}).
%
%
\begin{Thm}
\label{Thm:state}
\begin{enumerate}
\item
For  $(\ba,\bb)\in \Lambda(A)\times \Lambda(B)$,
$\ba\boxtimes \bb\in \Lambda(A\boxtimes B)$.
\item
Let $\rho_{\ba}$ be as in (\ref{eqn:kmsthree}).
For any $\ba,\bb\in \Lambda(*)$,
$\rho_{\ba}\ptimes \rho_{\bb}=\rho_{\ba\boxtimes \bb}$.
\item
Let $\varpi_{\ba}$ be as in Theorem \ref{Thm:typetwo}.
For any $\ba,\bb\in \Lambda(*)$,
$\varpi_{\ba}\ptimes \varpi_{\bb}$
is unitarily equivalent to $\varpi_{\ba\boxtimes \bb}$.
\end{enumerate}
\end{Thm}
From Theorem \ref{Thm:state}(i),
$\Lambda(*)$ is a semigroup with respect to the product $\boxtimes$.
By Theorem \ref{Thm:state}(ii),
the set $\{\rho_{\ba}:\ba\in\Lambda(*)\}$ of all states 
in (\ref{eqn:kmsthree}) is closed with respect to $\ptimes$,
and the mapping 
\[\Lambda(*)\ni\ba\mapsto \rho_{\ba}\in{\cal S}_{*}\] 
is a semigroup homomorphism from $(\Lambda(*),\boxtimes)$
to $({\cal S}_{*},\ptimes)$.
From Theorem \ref{Thm:state}(iii), we see that
%
%
\begin{equation}
\label{eqn:representation}
[\varpi_{\ba}]\ptimes [\varpi_{\bb}]=[\varpi_{\ba\boxtimes \bb}]
\quad(\ba,\bb\in\Lambda(*))
\end{equation}
where $[\pi]$ denotes the unitary equivalence class of $\pi\in \textsf{R}_{*}$.
From this, the mapping 
\[\Lambda(*)\ni\ba\mapsto [\varpi_{\ba}]\in{\cal R}_{*}\] 
is also a semigroup homomorphism from $(\Lambda(*),\boxtimes)$
to $({\cal R}_{*},\ptimes)$ for ${\cal R}_{*}$ in (\ref{eqn:rtwo}).

Next, we consider the type of $\varpi_{\ba}\ptimes \varpi_{\bb}$.
From Theorem \ref{Thm:state}(iii),
$(\varpi_{\ba}\ptimes \varpi_{\bb})(\co{A\boxtimes B})^{''}$ is 
also an AFD factor of type ${\rm III}_{\lambda}$
for $0<\lambda\leq 1$.
From Theorem \ref{Thm:typetwo} and \ref{Thm:state}(iii), the following holds.
%
%
\begin{cor}
\label{cor:tensor}
Let $\lambda(\ba)$ be as in (\ref{eqn:lambdatwo}).
Then
$(\varpi_{\ba}\ptimes \varpi_{\bb})(\co{A\boxtimes B})^{''}$ is 
an AFD factor of type ${\rm III}_{\lambda(\ba\boxtimes \bb)}$
for each $\ba,\bb\in\Lambda(*)$.
\end{cor}

\noindent
Remark that 
%
%
\begin{equation}
\label{eqn:lambda}
\lambda(\ba\boxtimes \bb)=\lambda(\bb\boxtimes \ba)
\quad(\ba,\bb\in\Lambda(*))
\end{equation}
by definition.

From Corollary \ref{cor:tensor}, 
the following holds.
%
%
\begin{cor}
\label{cor:type}
For $\be(A)$ in (\ref{eqn:eigentwo}),
$\varpi_{\be(A)}\ptimes \varpi_{\be(B)}$
is of type ${\rm III}_{\frac{1}{c_{A}c_{B}}}$
for each $A,B\in\cdm{*}$.
\end{cor}
%
%
\begin{rem}
\label{rem:tensor}
{\rm
\begin{enumerate}
\item
From Corollary \ref{cor:tensor},
we see that
the type of the tensor product $\varpi_{\ba}\ptimes \varpi_{\bb}$
is obtained by computing $\lambda(\ba\boxtimes \bb)$
for given $\ba,\bb$. However, it is not easy.
Remark that
$\lambda(\ba\boxtimes \bb)\ne
\lambda(\ba)\cdot \lambda(\bb)$ in general.
\item
From Lemma \ref{lem:factor},
we see that Theorem \ref{Thm:state}(iii), (\ref{eqn:lambda}) and
Corollary \ref{cor:type} are special cases
of tensor product of factor representations.
\item
We compare our results with those of tensor products
of type ${\rm III}_{\lambda}$ factors.
Let $N_{\lambda}$ denote
the AFD factors of type ${\rm III}_{\lambda}$ for $0<\lambda< 1$.
Then the following holds for any $\lambda$:
\[N_{\lambda}\otimes N_{\lambda}\cong N_{\lambda}. \]
We explain this in Appendix \ref{section:appone}.
On the other hand,
$\pi_{\lambda}\ptimes \pi_{\lambda}$
is of type ${\rm III}_{\lambda^{2}}$
when $\pi_{\lambda}=\varpi_{\be(A)}$ and 
$\lambda=1/c_{A}$
from Corollary \ref{cor:type}.
\item
The tensor product $\ptimes$ in (\ref{eqn:tensorstate}) was 
originally introduced 
such that it makes sense among permutative representations
of Cuntz-Krieger algebras \cite{TS05} but
it was not constructed as an operation among KMS states or 
type ${\rm III}$ representations.
However, Theorem \ref{Thm:state} sounds too good.
It glance at a certain system below.
We consider a relation between $\ptimes$ and $\boxtimes$ in 
Appendix \ref{section:apptwo}.
\item
Any KMS state in \cite{Okayasu} 
is written as $\rho_{\ba}$ in (\ref{eqn:kmsthree}) 
for a certain $\ba\in \Lambda(*)$.
This will be proved in $\S$ \ref{subsection:thirdfour}.
\item
In this paper, we treat special KMS states over Cuntz-Krieger algebras
but not general KMS states.
We do not know tensor products of general cases.
A slightly general case is considered in $\S$ \ref{subsection:thirdfive}.
\end{enumerate}
}
\end{rem}

In $\S$ \ref{section:second},
we review states and representations of C$^{*}$-algebras.
In $\S$ \ref{section:third},
we review KMS states over Cuntz-Krieger algebras
and prove main theorems in $\S$ \ref{subsection:firstfour}.
In $\S$ \ref{section:fourth}, we will compute tensor products
of type {\rm III} representations of Cuntz algebras. 
In $\S$ \ref{section:fifth}, we will show formulae of tensor powers
of representations.

%
%
\sftt{States and representations of C$^{*}$-algebras}
\label{section:second}
In this section,
we review states and representations of C$^{*}$-algebras.

%
%
\ssft{Equivalences and the type of a representation}
\label{subsection:secondone}
We review factor representations, 
two equivalences of representations 
and the type of a representation \cite{Blackadar2006,Dixmier,Ped,Wallach}.

Let ${\cal A}$ be a C$^{*}$-algebra.
A nondegenerate representation $({\cal H},\pi)$ of ${\cal A}$ 
is called a {\it factor (or primary or factorial) representation} 
if $\pi({\cal A})^{''}$ is a factor.
For two representations $\pi_{1}$ and $\pi_{2}$ of ${\cal A}$,
we state that $\pi_{1}$ and $\pi_{2}$ are {\it unitarily equivalent 
(or spatially equivalent)}
if there exists a unitary $u$
such that $\pi_{2}(x)=u\pi_{1}(x)u^{*}$ for each $x\in {\cal A}$;
$\pi_{1}$ and $\pi_{2}$ are {\it quasi-equivalent} 
if there exists 
a  $*$-isomorphism $f$ from $\pi_{1}({\cal A})^{''}$ 
to $\pi_{2}({\cal A})^{''}$
such that $f(\pi_{1}(x))=\pi_{2}(x)$ for each $x\in {\cal A}$;
$\pi_{1}$ and $\pi_{2}$ are {\it disjoint} if no subrepresentation
of $\pi_{1}$ is equivalent to a subrepresentation of $\pi_{2}$.

Let $\simeq$ and $\approx$ be as in $\S$ \ref{subsection:firstone}.
We review known results as follows:
For two representations $\pi_{1}$ and $\pi_{2}$,
$\pi_{1}\approx \pi_{2}$
if and only if
$\pi_{1}\simeq \pi_{2}$ up to multiplicity.
Especially,
when $\pi_{1}$ and $\pi_{2}$ are irreducible,
$\pi_{1}\approx \pi_{2}$
if and only $\pi_{1}\simeq \pi_{2}$.
Two factor representations are either disjoint or quasi-equivalent. 
Any irreducible representation is a factor representation.
For a set $\{\pi_{i}\}$ of pairwise disjoint representations 
of a C$^{*}$-algebra ${\cal A}$,
$(\oplus\pi_{i})({\cal A})^{''}=\oplus\pi_{i}({\cal A})^{''}$.

A nondegenerate representation $\pi$ of a C$^{*}$-algebra ${\cal A}$
is said to be {\it (pure) type {\rm X}}
if $\pi({\cal A})^{''}$ is of type X 
for X$=$I, II, III, II$_{1}$, II$_{\infty}$.
The type of a representation is invariant under quasi-equivalence.
%
%
\begin{lem}
\label{lem:three}
(\cite{Dixmier}, $\S$ 5.6.6)
Let $\pi_{1},\pi_{2}$ be two representations of a C$^{*}$-algebra
in separable spaces and assume that 
$\pi_{1}$ or $\pi_{2}$ is of type ${\rm III}$.
Then $\pi_{1}\approx \pi_{2}$ implies $\pi_{1}\simeq\pi_{2}$.
\end{lem}

%
%
\ssft{KMS states}
\label{subsection:secondtwo}
In this subsection, 
we review basic facts of KMS states according to \cite{BR2}.
%
%
\begin{defi}
\label{defi:kmsone}
\begin{enumerate}
\item
Let $({\cal A}, {\bf R},\tau)$ be a C$^{*}$-dynamical system. 
The state $\rho$ over ${\cal A}$ is a {\it $\tau$-KMS 
state at value $\beta\in {\bf R}$, or a $(\tau, \beta)$-KMS state, }
if $\rho(a\tau_{i\beta}(b)) = \rho(ba)$
for all $a,b$ in a norm dense, $\tau$-invariant $*$-subalgebra of ${\cal A}$.
We call $\beta$ {\it the inverse temperature of $\rho$.}
\item
A state $\rho$ over a C$^{*}$-algebra ${\cal A}$
is a {\it KMS state} if there exists an action $\tau$ of ${\bf R}$ 
on ${\cal A}$ and $\beta\in {\bf R}$ such that $\rho$ 
is a $(\tau,\beta)$-KMS state.
\end{enumerate}
\end{defi}

\noindent
Note that a $(\tau_{\theta t}, \beta)$-KMS state coincides with
a $(\tau_{\theta^{'}t},\beta^{'})$-KMS state 
when $\theta\beta=\theta^{'}\beta^{'}$.
From this,
$(\tau_{\gamma t}, \gamma^{-1}\beta)$-KMS state $\rho$ 
is independent in the choice of 
$\gamma\in {\bf R}\setminus\{0\}$.

A state $\rho$ is a {\it factor state (or factorial state
or primary state)}
if the GNS representation of $\rho$ is a factor representation
\cite{Blackadar2006,Dixmier}.
Let $K_{\beta}(\tau)$ denote the set of all $(\tau,\beta)$-KMS states.
For $\rho\in K_{\beta}(\tau)$,
$\rho$ is an extremal point of $K_{\beta}(\tau)$
if and only if $\rho$ is a factor state.
From this, if $(\tau,\beta)$-KMS state exists uniquely,
then it is a factor state.

We prepare a lemma about C$^{*}$-subalgebras and their representations. 
%
%
\begin{lem}
\label{lem:kmscor}
(\cite{Izumi}, Remark 4.2)
Let $\rho$ be a KMS state over a unital C$^{*}$-algebra ${\cal B}$
and let ${\cal C}$ be a unital C$^{*}$-subalgebra of ${\cal B}$.
Let $\pi_{\rho}$ and $\pi_{\rho|_{{\cal C}}}$
denote the GNS representations of ${\cal B}$ and ${\cal C}$ with respect to
$\rho$ and $\rho|_{{\cal C}}$, respectively.
Then $(\pi_{\rho})|_{{\cal C}}$ and $\pi_{\rho|_{{\cal C}}}$ 
are quasi-equivalent.
\end{lem}

%
%
\sftt{KMS states over Cuntz-Krieger algebras}
\label{section:third}
We review KMS states and their GNS representations of Cuntz-Krieger algebras
in this section.
Proofs of main theorems in $\S$ \ref{subsection:firstfour} 
will be given in $\S$ \ref{subsection:thirdfour}.
Let $\ndm{n}$ and $\cdm{*}$
be as in $\S$ \ref{subsection:firstone} and (\ref{eqn:nondegenerate}), respectively.

%
%
\ssft{Perron-Frobenius theorem}
\label{subsection:thirdone}
If a non-negative matrix $A$ is irreducible, 
the Perron-Frobenius theorem guarantees 
the existence of the strictly positive eigenvector with
respect to the simple root $a$ of the characteristic polynomial such that 
$a\geq |b|$ for any other eigenvalue $b$.
From the Perron-Frobenius theorem,
the statement at the beginning of $\S$ \ref{subsection:firstthree}
is rewritten as follows.
%
%
\begin{lem}
\label{lem:eigenfour}
For any $A\in\cdm{n}$ and $\bomega=(\omega_{1},\ldots,\omega_{n}),
\in {\bf R}_{+}^{n}$,
there exists  unique $\beta>0$ such that the vector
%
%
\begin{equation}
\label{eqn:relation}
\ba\equiv (e^{-\beta\omega_{1}},\ldots,e^{-\beta\omega_{n}})
\end{equation}
belongs to $\Lambda(A)$ in (\ref{eqn:pfe}).
\end{lem}

%
%
\ssft{Cuntz-Krieger algebras}
\label{subsection:thirdtwo}
For $A=(A_{ij})\in \ndm{n}$,
$\coa$ is the
{\it \ck\ algebra by $A$} if 
$\coa$ is a C$^{*}$-algebra
which is universally generated by
partial isometries $s_{1},\ldots,s_{n}$
and they satisfy $s_{i}^{*}s_{i}=\sum_{j=1}^{n}A_{ij}s_{j}s_{j}^{*}$
for $i=1,\ldots,n$ and $\sum_{i=1}^{n}s_{i}s_{i}^{*}=I$ \cite{CK}.
The C$^{*}$-algebra $\coa$ is simple 
if and only if $A\in\cdm{*}$.
For $\boxtimes$ in (\ref{eqn:phitwo}),
if both $\coa$ and $\co{B}$ are simple, then so is $\co{A\boxtimes B}$.

For methods to construct representations of $\coa$,
see \cite{CKR01,CKR02,PFO01}.
For type ${\rm I}$ representations (especially, irreducible representations) of $\coa$,
see \cite{CKR02}. 
There exists no type ${\rm II}$ representation
of $\coa$ when $\coa$ is simple
because $\coa$ is purely infinite (\cite{RS}, R\o rdam, Proposition 4.4.2) and 
a purely infinite C$^{*}$-algebra has no nondegenerate lower semicontinuous trace
(\cite{Blackadar2006}, Proposition V.2.2.29).
For embeddings of Cuntz-Krieger algebras, see \cite{CK01}.

%
%
\ssft{KMS states over Cuntz-Krieger algebras}
\label{subsection:thirdthree}
According to \cite{EL2,Okayasu},
we review certain KMS states and their representations of Cuntz-Krieger algebras.
For $\bomega=(\omega_{1},\ldots,\omega_{n})\in {\bf R}^{n}_{+}$
and $t\in {\bf R}$,
define the $*$-automorphism $\alpha_{t}^{\bomega}$ of $\coa$ by
%
%
\begin{equation}
\label{eqn:automorphism}
\alpha_{t}^{\bomega}(s_{i})
\equiv e^{\sqrt{-1}\omega_{i}t}s_{i}\quad(i=1,\ldots,n).
\end{equation}
Then $\alpha^{\bomega}$ is a one-parameter automorphism group of $\coa$.
%
%
\begin{Thm}
\label{Thm:exellaca}
(\cite{EL2}, Proposition 18.3, Theorem 18.5)
Assume that $A\in\cdm{n}$ and $\alpha^{\bomega}$ 
is as in (\ref{eqn:automorphism}).
Choose $\beta>0$ for $A$ and $\bomega$ is as in Lemma \ref{lem:eigenfour}.
Then an $\alpha^{\bomega}$-KMS state $\phi^{\bomega}$ 
over $\coa$ with an inverse temperature $\beta$
is given as follows:
%
%
\begin{equation}
\label{eqn:kmsone}
\phi^{\bomega}(s_{J}s_{K}^{*})=
\delta_{JK}e^{-\beta\omega_{j_{1}}}\cdots
e^{-\beta\omega_{j_{m-1}}}x_{j_{m}}
\end{equation}
when $s_{J}s_{K}^{*}\ne 0$
for $J=(j_{1},\ldots,j_{m})\in\{1,\ldots,n\}^{m}$
and $K\in\cup_{l\geq 1}\{1,\ldots,n\}^{l}$
where $\bx=(x_{1},\ldots,x_{n})$ is the Perron-Frobenius eigenvector of 
$\hat{\ba}A$ for $\ba$ in  Lemma \ref{lem:eigenfour}
such that $x_{1}+\cdots+x_{n}=1$.
Furthermore,
an $\alpha^{\bomega}$-KMS state is unique and $\beta$ is also unique.
\end{Thm}
%
%
\begin{Thm}
\label{Thm:oka}
(\cite{Okayasu}, Theorem 4.2)
In addition to assumptions in Theorem \ref{Thm:exellaca},
let $\pi_{\phi^{\bomega}}$ denote the GNS representation 
of $\coa$ by $\phi^{\bomega}$,
$M\equiv \pi_{\phi^{\bomega}}(\coa)^{''}$ and
let $G$ denote the closed subgroup of the additive group ${\bf R}$
generated by $\beta\omega_{i}$ for all $i$.
\begin{enumerate}
\item
If $\omega_{i}/\omega_{j}\in {\bf Q}$ for all $i,j\in \{1,\ldots,n\}$,
then $M$ is an AFD factor of type {\rm III}$_{\lambda}$ for $\lambda=e^{-r}$,
where $G=r{\bf Z}$ for some $r\in {\bf R}_{+}$.
\item
If $\omega_{i}/\omega_{j}\not\in {\bf Q}$ 
for some $i,j\in\{1,\ldots,n\}$,
then $M$ is an AFD factor of type {\rm III}$_{1}$.
\end{enumerate}
\end{Thm}

%
%
\ssft{Proofs of main theorems}
\label{subsection:thirdfour}
We prove main theorems in this subsection.\\

\noindent
{\it Proof of Lemma \ref{lem:factor}.}
(i)
Let $\pi_{12}$ be a representation of $\co{2}$
with a cyclic vector $\Omega$ such that
$\pi_{12}(s_{1}s_{2})\Omega=\Omega$.
Then such a representation exists uniquely up to unitary equivalence
and it is irreducible.
Let $P_{2}(12)$ denote the unitary equivalence class of $\pi_{12}$.
By $\S$ 4.1 in \cite{TS01},
the irreducible decomposition 
of $k$ times tensor power $P_{2}(12)^{\ptimes k}$ of $P_{2}(12)$ 
is multiplicity-free with $2^{k-1}$ irreducible components for each $k\geq 1$.
Since $P_{2}(12)$ is irreducible, it is a (class of) factor representation.
However, $P_{2}(12)^{\ptimes k}$ is not when $k\geq 2$.

\noindent 
(ii)
For $i=1,\ldots,n$,
let $\pi_{i}$ be a representation of $\con$
with a cyclic vector $\Omega_{i}$ such that
$\pi_{i}(s_{i})\Omega_{i}=\Omega_{i}$.
Then such a representation exists uniquely up to unitary equivalence
and it is irreducible.
Furthermore $\pi_{i}\not\simeq\pi_{j}$ when $i\ne j$.
Let $P_{n}(i)$ denote the unitary equivalence class of $\pi_{i}$.
By $\S$ 4.1 in \cite{TS01}, 
$P_{2}(1)\ptimes P_{2}(2)=P_{4}(2)$ and 
$P_{2}(2)\ptimes P_{2}(1)=P_{4}(3)$.
Since both $P_{4}(2)$ and $P_{4}(3)$ are
irreducible and $P_{4}(2)\not\simeq P_{4}(3)$,
$P_{4}(2)\not\approx P_{4}(3)$.
Hence the statement holds.

\noindent 
(iii)
By the definition of $\ptimes$,
the statement holds.

\noindent 
(iv)
Let $\pi_{12}$ and $\Omega$ be as in the proof of (i).
Assume $\|\Omega\|=1$.
Define the state $\rho$ over $\co{2}$ by
$\rho\equiv \langle\Omega|\pi_{12}(\cdot)\Omega \rangle$.
Then we see that 
\[\rho\ptimes \rho=
\langle\Omega\otimes \Omega|
(\pi_{12}\ptimes \pi_{12})(\cdot)\Omega\otimes \Omega \rangle.\]
On the other hand,
$\pi_{12}\ptimes \pi_{12}$ has just two irreducible components
and is multiplicity-free, and $\pi_{12}\simeq\pi_{\rho}$.
On the other hand,
$\pi_{\rho\ptimes \rho}$ is irreducible.
Hence the statement holds.
\qedh

\noindent
Remark that the statement in $\S$ 3.2 of \cite{TS01} is wrong.
In general,
$(\pi_{\rho_{1}}\ptimes \pi_{\rho_{2}})|_{{\cal K}}$
and $\pi_{\rho_{1}\ptimes \rho_{2}}$
are unitarily equivalent
where ${\cal K}$ denotes the cyclic representation space
with the tensor product $\Omega_{\rho_{1}}\otimes \Omega_{\rho_{2}}$
of GNS cyclic vectors as the cyclic vector. 
\\

\noindent
{\it Proof of Theorem \ref{Thm:state}.}
(i)
For any $A\in M_{n}({\bf C}), B\in M_{m}({\bf C})$
and $\bx\in {\bf C}^{n},\by\in {\bf C}^{m}$,
we see that
$\widehat{\bx\boxtimes \by}(A\boxtimes B)
=\hat{\bx}A\boxtimes \hat{\by}B\in M_{nm}({\bf C})$
and $(A\boxtimes B)(\bx\boxtimes \by)=A\bx \boxtimes B\by\in {\bf C}^{nm}$.
From these,
if $\bx$ and $\by$ are the Perron-Frobenius eigenvectors of 
$\hat{\ba}A$ and $\hat{\bb}B$, respectively, then
$\bx\boxtimes \by$ is also that of the matrix $\hat{\ba}A\boxtimes \hat{\bb}B$.
From this, the statement holds.

\noindent
(ii)
From (i), $\rho_{\ba\boxtimes \bb}$ is well-defined.
By definitions of $\ptimes$ and $\rho_{\ba}$,
the statement is verified directly.

\noindent
(iii)
For $\ba= (a_{1},\ldots,a_{n})\in\Lambda(A)$,
define
$\bomega\equiv (-\log a_{1},\ldots,-\log a_{n})$.
Then $\rho_{\ba}$ in (\ref{eqn:kmsthree})
coincides with $\phi^{\bomega}$ in (\ref{eqn:kmsone}) 
and the inverse temperature $1$.
Therefore $\rho_{\ba}$ is the unique 
$(\alpha^{\bomega},1)$-KMS state over $\coa$.

Assume that $\ba\in \Lambda(A)$ and $\bb\in\Lambda(B)$.
Let ${\cal C}\equiv \varphi_{A,B}(\co{A\boxtimes B})$,
${\cal B}\equiv 
\co{A}\otimes \co{B}$ and $\rho\equiv \rho_{\ba}\otimes \rho_{\bb}$.
Then $\rho$ is a state over the unital C$^{*}$-algebra ${\cal B}$.
By the definition of $\ptimes$,
$\rho_{\ba}\ptimes \rho_{\bb}=\rho|_{{\cal C}}$.
From (ii), $\rho_{\ba\boxtimes \bb}=\rho|_{{\cal C}}$.
By this and the definition of $\varpi_{\ba}$,
\[\varpi_{\ba\boxtimes \bb}=\pi_{\rho|_{{\cal C}}}.\]
On the other hand,
$\varpi_{\ba}\ptimes \varpi_{\bb}
=(\varpi_{\ba}\otimes \varpi_{\bb})|_{{\cal C}}$.
Since $\varpi_{\ba}\otimes \varpi_{\bb}$ is unitarily equivalent to
the GNS representation $\pi_{\rho}$ of ${\cal B}$ by $\rho$,
we can identify $\pi_{\rho}$ with $\varpi_{\ba}\otimes \varpi_{\bb}$.
Hence
\[\varpi_{\ba}\ptimes \varpi_{\bb}=(\pi_{\rho})|_{{\cal C}}.\]
From Lemma \ref{lem:kmscor}, 
$\varpi_{\ba}\ptimes \varpi_{\bb}\approx \varpi_{\ba\boxtimes \bb}$.
By Theorem \ref{Thm:typetwo} and Lemma \ref{lem:three},
the statement holds.
\qedh

In Theorem \ref{Thm:exellaca},
let $\ba$ be as in (\ref{eqn:relation}).
Then $\ba$ belongs to $\Lambda(A)$ in (\ref{eqn:pfe}).
We see that $\rho_{\ba}$ in (\ref{eqn:kmsthree})
coincides with $\phi^{\bomega}$ in (\ref{eqn:kmsone}).
From this and the proof of Theorem \ref{Thm:state}(ii),
the set of all KMS states in Theorem \ref{Thm:oka} coincides with
the set $\{\rho_{\ba}:\ba\in\Lambda(*)\}$.

\ww
{\it Proof of Corollary \ref{cor:type}.}
Since $c_{A\boxtimes B}=c_{A}c_{B}$,
$\be(A\boxtimes B)=\be(A)\boxtimes \be(B)$.
This implies the statement. 
\qedh

%
%
\ssft{Tensor products of general KMS states}
\label{subsection:thirdfive}
We consider the sufficient condition such that
the tensor product of two (slightly) general KMS states is also a KMS state.
%
%
\begin{lem}
\label{lem:kmsproduct}
Let $(\co{A_{i}},{\bf R},\alpha^{(i)})$ be
a C$^{*}$-dynamical system 
and let $\rho_{i}$ be an $\alpha^{(i)}$-KMS state over $\co{A_{i}}$ 
with an inverse temperature $\beta_{i}\ne 0$ such that 
the KMS condition holds on an $\alpha^{(i)}$-invariant dense $*$-subalgebra ${\cal B}_{i}$
of $\co{A_{i}}$ for $i=1,2$.
We assume the following conditions:
\begin{enumerate}
\item
there exists a dense $*$-subalgebra ${\cal C}$
of $\co{A_{1}\boxtimes A_{2}}$ such that 
$\varphi_{A_{1},A_{2}}({\cal C})$ is included in the
algebraic tensor product ${\cal B}_{1}\odot {\cal B}_{2}$
where $\varphi_{A_{1},A_{2}}$ is as in (\ref{eqn:varphi}),
\item
$(\alpha^{(1)}\otimes \alpha^{(2)})(\varphi_{A_{1},A_{2}}({\cal C}))
\subset \varphi_{A_{1},A_{2}}({\cal C})$.
\end{enumerate}
Define the action $\alpha^{(1)}\ptimes \alpha^{(2)}$ of ${\bf R}$
on $\co{A_{1}\otimes A_{2}}$ by 
\[
(\alpha^{(1)}\ptimes \alpha^{(2)})_{t}
\equiv \varphi_{A_{1},A_{2}}^{-1}
\circ (\alpha^{(1)}_{\beta_{1}t}\otimes \alpha^{(2)}_{\beta_{2}t})
\circ \varphi_{A_{1},A_{2}}\quad(t\in {\bf R}).
\]
Then $\rho_{1}\ptimes \rho_{2}$
is an $\alpha^{(1)}\ptimes \alpha^{(2)}$-KMS state 
over $\co{A_{1}\boxtimes A_{2}}$ with
an inverse temperature $1$.
\end{lem}
%
%
\pr
Let $x,y\in {\cal C}$.
By assumption,
there exist 
$x_{1}^{'},\ldots,x_{l}^{'},y_{1}^{'},\ldots,y_{k}^{'}\in {\cal B}_{1}$
and $x_{1}^{''},\ldots,x_{l}^{''},
y_{1}^{''},\ldots,y_{k}^{''}\in {\cal B}_{2}$
such that
$\varphi_{A_{1},A_{2}}(x)=\sum_{i}x_{i}^{'}\otimes x_{i}^{''}$
and $\varphi_{A_{1},A_{2}}(y)=\sum_{i}y_{i}^{'}\otimes y_{i}^{''}$.
Then
we can verify that
%
%
\begin{equation}
\label{eqn:auto}
(\rho_{1}\ptimes \rho_{2})(x\,(\alpha^{(1)}\ptimes \alpha^{(2)})_{t}(y))
=\sum_{i,j}
\rho_{1}(x_{i}^{'}\alpha_{\beta_{1}t}^{(1)}(y_{j}^{'}))
\,\rho_{2}(x_{i}^{''}\alpha_{\beta_{2}t}^{(2)}(y_{j}^{''}))
\end{equation}
for each $t\in {\bf R}$.
On the other hand,
%
%
\begin{equation}
\label{eqn:autotwo}
\sum_{i,j}
\rho_{1}(y_{j}^{'}x_{i}^{'})
\,\rho_{2}(y_{j}^{''}x_{i}^{''})=(\rho_{1}\ptimes \rho_{2})(yx).
\end{equation}
We obtain the KMS condition for $\rho_{1}\ptimes \rho_{2}$
on ${\cal C}$ by taking the limit $t\to \sqrt{-1}$ in (\ref{eqn:auto}).
\qedh

\noindent
From the definition of $\varphi_{A,B}$ in (\ref{eqn:varphi}) and 
Lemma \ref{lem:kmsproduct},
the following holds.
%
%
\begin{cor}
\label{cor:product}
In addition to Lemma \ref{lem:kmsproduct},
if ${\cal B}_{1}$ and ${\cal B}_{2}$
are $*$-subalgebras generated by
canonical generators of $\co{A_{1}}$
and $\co{A_{2}}$, respectively,
then ${\cal C}$ can be taken as the 
dense $*$-subalgebra generated by
canonical generators of $\co{A_{1}\boxtimes A_{2}}$.
\end{cor}

%
%
\sftt{Cases of Cuntz algebras}
\label{section:fourth}
We show cases of Cuntz algebras in this section.
%
%
\ssft{KMS states over $\con$}
\label{subsection:fourthone}
GNS representations of 
KMS states over Cuntz algebras were studied by \cite{Izumi}.
We rewrite them as special cases of Cuntz-Krieger algebras.
For $n\geq 2$,
let ${\rm Int}\Delta_{n-1}$ denote the interior of the $n-1$-simplex, that is,
${\rm Int}\Delta_{n-1}\equiv
\{(a_{1},\ldots,a_{n})\in {\bf R}^{n}
:\sum_{j=1}^{n}a_{j}=1,\,a_{i}>0\mbox{ {\rm for all} }i\}$.
Define $F_{n}\in\cdm{n}$ by 
%
%
\begin{equation}
\label{eqn:fn}
(F_{n})_{ij}=1\quad(i,j=1,\ldots,n).
\end{equation}
Let $\Lambda(A)$ be as in (\ref{eqn:pfe}).
Then we see that $\Lambda(F_{n})\subset {\rm Int}\Delta_{n-1}$.
Since an eigenvalue of $\hat{\ba}F_{n}$ is $0$ or $1$ 
for each $\ba\in {\rm Int}\Delta_{n-1}$,
the following holds.
%
%
\begin{lem}
\label{lem:lambdaf}
(\cite{Izumi}, $\S$ 2.1)
For each $n\geq 2$, $\Lambda(F_{n})={\rm Int}\Delta_{n-1}$.
\end{lem}

For a finite set $\{p_{1},\ldots,p_{n}\}$ of natural numbers,
let $\gcd\{p_{1},\ldots,p_{n}\}$ denote
the greatest common divisor of $\{p_{1},\ldots,p_{n}\}$. 
From Theorem \ref{Thm:typetwo} and Lemma \ref{lem:lambdaf},
the following holds.
%
%
\begin{prop}
\label{prop:cuntzcase}
If $p_{1},\ldots,p_{n}\in {\bf N}$ 
and $\lambda >0$ satisfy that
$\gcd\{p_{1},\ldots,p_{n}\}=1$ and 
\[\lambda^{p_{1}}+\cdots+\lambda^{p_{n}}=1,\]
then the representation
$\varpi_{(\lambda^{p_{1}},\ldots,\lambda^{p_{n}})}$ of $\con$ 
is of type ${\rm III}_{\lambda}$.
\end{prop}

\noindent
From Proposition \ref{prop:cuntzcase},
$\lambda$ is an algebraic number
when $\varpi_{\ba}$ is of type ${\rm III}_{\lambda}$.
This was pointed out by Theorem 4.7(ii) in \cite{Izumi}.

For $\ba\in \Lambda(F_{n})$,
$\ba$ is also the Perron-Frobenius eigenvector of $\hat{\ba}F_{n}$
which satisfies the assumption for $\bx$ in (\ref{eqn:kmsthree}).
Hence the KMS state $\rho_{\ba}$ over $\con$ in (\ref{eqn:kmsthree})
is given as follows:
%
%
\begin{equation}
\label{eqn:kmsfour}
\rho_{\ba}(s_{J}s_{K}^{*})=\delta_{JK}a_{J}
\quad(J,K\in \cup_{l\geq 0}\{1,\ldots,n\}^{l})
\end{equation}
where $s_{1},\ldots,s_{n}$ denote
canonical generators of $\con$ 
and $a_{J}\equiv a_{j_{1}}\cdots a_{j_{m}}$
when $J=(j_{1},\ldots,j_{m})$.
This type state is called {\it quasi-free} \cite{Evans}.
Especially,
when $\ba=(\frac{1}{n},\ldots,\frac{1}{n})$, 
we write $\rho^{(n)}$ as $\rho_{\ba}$. 
Then the following holds:
\begin{enumerate}
\item
(\cite{BR2}, Example 5.3.27)
$\rho^{(n)}(s_{J}s_{K}^{*})=\delta_{J,K}n^{-|J|}$ for $J,K\in\cup_{l\geq 0}\nset{l}$
where $|J|\equiv k$ for $J=(j_{1},\ldots,j_{k})$.
\item
$\rho^{(n)}\ptimes \rho^{(m)}=\rho^{(nm)}$
for each $n,m\geq 2$.
\end{enumerate}

%
%
\ssft{Formulae of tensor product}
\label{subsection:fourthtwo}
We show several formulae of tensor products 
of (unitary equivalence classes of) representations 
of Cuntz algebras in this subsection.
From Corollary \ref{cor:tensor}, the essential computation is reduced to 
that of $\lambda(\ba\boxtimes \bb)$.
For $F_{n}$ in (\ref{eqn:fn}),
we compute 
$\lambda(\ba\boxtimes \bb)$ for $\ba,\bb\in\Lambda(F_{n})$.

%
%
\begin{ex}
\label{ex:formulae}
{\rm
We consider the case of $\co{2}$.
Define $\ba,\bb,\bc\in \Lambda(F_{2})$ by
\[\ba=\left(\frac{1}{3},\frac{2}{3}\right),\quad
\bb=\left(\frac{1}{2},\frac{1}{2}\right),\quad
\bc=\Bigl(\frac{\sqrt{5}-1}{2},\Bigl\{\frac{\sqrt{5}-1}{2}\Bigr\}^{2}\Bigr).
\]
From Theorem \ref{Thm:typetwo} and Lemma \ref{lem:log},
$\lambda(\ba)=1$,
$\lambda(\bb)=\frac{1}{2}$ and 
$\lambda(\bc)=\frac{\sqrt{5}-1}{2}$.
Hence 
$\varpi_{\ba}$, $\varpi_{\bb}$ and $\varpi_{\bc}$ 
are of type ${\rm III}_{1}$, ${\rm III}_{\frac{1}{2}}$,${\rm III}_{\frac{\sqrt{5}-1}{2}}$,
respectively.
\begin{enumerate}
\item
Since $\ba\boxtimes \bb
=(\frac{1}{6},\frac{1}{6},\frac{1}{3},\frac{1}{3})$
and $\log \frac{1}{6}/\log \frac{1}{3}\not\in {\bf Q}$,
$\varpi_{\ba}\ptimes \varpi_{\bb}$ is of type ${\rm III}_{1}$. 
\item
Since
$\bb\boxtimes\bc=(\frac{\sqrt{5}-1}{4},
\frac{(\sqrt{5}-1)^{2}}{8},\frac{\sqrt{5}-1}{4},
\frac{(\sqrt{5}-1)^{2}}{8})$ and
 $\log \frac{\sqrt{5}-1}{4}/\log(\frac{(\sqrt{5}-1)^{2}}{8})
\not\in {\bf Q}$,
$\varpi_{\bb}\ptimes \varpi_{\bc}$
is of type ${\rm III}_{1}$.
\end{enumerate}
}
\end{ex}

\noindent
Furthermore,
we show the following.
%
%
\begin{prop}
\label{prop:oneone}
For any $n,m\geq 2$,
there exists $\{\ba_{n}:n\geq 2\}$
such that $\ba_{n}\in\Lambda(F_{n})$,
both $\varpi_{\ba_{n}}$ and
$\varpi_{\ba_{n}}\ptimes \varpi_{\ba_{m}}$ are of type ${\rm III}_{1}$
for each $n,m\geq 2$.
\end{prop}
%
%
\pr
We prepare the following fact:
%
%
\begin{equation}
\label{eqn:abs}
\mbox{If $k$ is an odd integer greater than equal $3$, then
$\log k/\log \frac{k}{2}\not\in {\bf Q}$.}
\end{equation}
This can be proved by reduction to absurdity.

For $n\geq 2$,
define $\ba_{n}\in\Lambda(F_{n})$ by
\[
\ba_{n}\equiv \left\{
\begin{array}{ll}
\Bigl(\underbrace{\frac{1}{n+1},\ldots,
\frac{1}{n+1}}_{n-1\mbox{ {\footnotesize times}}},
\disp{\frac{2}{n+1}}\Bigr)
\quad &
\mbox{when $n$ is even},\\
\\
\Bigl(\underbrace{\frac{1}{n+2},\ldots,
\frac{1}{n+2}}_{n-2\mbox{ {\footnotesize times}}},
\disp{\frac{2}{n+2},
\frac{2}{n+2}}\Bigr)\quad &
\mbox{when $n$ is odd.}
\end{array}
\right.
\]
From Theorem \ref{Thm:typetwo}(ii), Lemma \ref{lem:log} and (\ref{eqn:abs}),
$\varpi_{\ba_{n}}$ is of type ${\rm III}_{1}$ for each $n\geq 2$.
By definition,
there exists an odd integer $k$ which is greater than equal $9$, and 
any component of $\ba_{n}\boxtimes \ba_{m}$ is 
$\frac{1}{k}$ or $\frac{2}{k}$ or $\frac{4}{k}$.
From Lemma \ref{lem:log} and (\ref{eqn:abs}),
$\lambda(\ba_{n}\boxtimes \ba_{m})=1$.
Since 
$\varpi_{\ba_{n}}\ptimes \varpi_{\ba_{m}}$ is of type 
${\rm III}_{\lambda(\ba_{n}\boxtimes \ba_{m})}$,
the statement holds.
\qedh

%
%
\sftt{Tensor powers of representations}
\label{section:fifth}
We consider the power of tensor product (= tensor power)
of representations in this section.
For a (unitary equivalence class of) representation $\pi$ of $\coa$,
let $\pi^{\ptimes k}$ denote $\pi\ptimes \cdots\ptimes \pi$ ($k$ times).
For $\ba\in \Lambda(A)$,
let $\lambda(\ba)$ be as in (\ref{eqn:lambdatwo}).
From Corollary \ref{cor:tensor},
$\varpi_{\ba}^{\ptimes k}$ is of type 
${\rm III}_{\lambda(\ba^{\boxtimes k})}$ for $k\geq 1$
where $\ba^{\boxtimes k}\equiv \ba\boxtimes \cdots\boxtimes \ba$ ($k$ times).
Remark that $\lambda(\ba^{\boxtimes k})\ne\lambda(\ba)^{k}$ in general.
From Corollary \ref{cor:type},
$\varpi_{\be(A)}^{\ptimes k}$
is of type ${\rm III}_{(1/c_{A})^{k}}$ for each $k\geq 1$.

%
%
\ssft{Tensor powers of type ${\rm III}_{1}$ factor representations}
\label{subsection:fifthone}
%
%
\begin{prop}
\label{prop:onetwo}
For any $A\in\cdm{*}$ and $\ba\in\Lambda(A)$,
if $\varpi_{\ba}$ is of type ${\rm III}_{1}$,
then
$\varpi_{\ba}^{\ptimes k}$ is also of type ${\rm III}_{1}$
for each $k\geq 1$.
\end{prop}
%
%
\pr
Assume $\ba=(a_{1},\ldots,a_{n})$.
By assumption and Lemma \ref{lem:log},
there exist $i,j$ such that 
$\log a_{i}/\log a_{j}\not\in {\bf Q}$.
On the other hand,
$a_{i}^{k}$ and $a_{j}^{k}$ always appear as
components of $\ba^{\boxtimes k}$.
Then $\log a_{i}^{k}/\log a_{j}^{k}=\log a_{i}/\log a_{j}\not\in {\bf Q}$.
From this, 
$\lambda(\ba^{\boxtimes k})=1$ for any $k\geq 1$.
Hence the statement holds.
\qedh

%
%
\ssft{Periodicity of the tensor power of a representation of $\co{2}$}
\label{subsection:fifthtwo}
We show formulae of tensor power of representations of $\co{2}$
in this subsection.
For $\Lambda(*)$ in (\ref{eqn:pfe}) and
$F_{n}$ in (\ref{eqn:fn}),
let $\Lambda_{2}\equiv \Lambda(F_{2})$.
From Lemma \ref{lem:lambdaf},
$\Lambda_{2}=\{(x,1-x)\in {\bf R}_{+}^{2}: 0< x<1\}$.
For $\ba\in\Lambda_{2}$,
let $\rho_{\ba}$ be as in (\ref{eqn:kmsfour}).
For any $(a,b)\in \Lambda_{2}$,
let $\varpi_{(a,b)}$ denote the GNS representation of $\co{2}$ by $\rho_{(a,b)}$.
We compute the type of $\varpi_{(a,b)}^{\ptimes k}$.
%
%
\begin{prop}
\label{prop:power}
If $x>0$ and $p,q\in {\bf N}$ satisfy that 
\[\gcd\{p,q\}=1\quad\mbox{ and }\quad x^{p}+x^{q}=1,\]
 then
$\varpi_{(x^{p},x^{q})}^{\ptimes k}$
is of {\rm type} ${\rm III}_{x^{r}}$ for each $k\geq 1$
where
\[r\equiv \gcd\{|p-q|,k\}\]
and we define $\gcd\{k,0\}\equiv k$ for convenience.
\end{prop}
%
%
\pr
Since $\varpi_{(a,b)}^{\ptimes k}$ is of type 
${\rm III}_{\lambda((a,b)^{\boxtimes k})}$,
we compute $\lambda((a,b)^{\boxtimes k})$ as follows.
Let $\ba=(x^{p},x^{q})$.
Assume $p=q$.
Since $\gcd\{p,q\}=1$, $(p,q)=(1,1)$.
Then $\ba^{\boxtimes k}=(x^{k},\ldots,x^{k})$.
Therefore $\lambda(\ba^{\boxtimes k})=x^{k}$.
Hence the statement holds.

Assume $t\equiv p-q>0$
and $r\equiv \gcd\{k,t\}$.
By assumption, $\gcd\{t,q\}=1$.
Then we can write $k=rk_{0}$ and $t=rt_{0}$
for some $k_{0}, t_{0}\in {\bf N}$
such that $\gcd\{k_{0},t_{0}\}=1$.
For any component of $\ba^{\boxtimes k}$,
there exists $i\in\{0,\ldots,k\}$ such that it is written as  
$x^{pi+q(k-i)}$.
Define $l_{i}\equiv pi+q(k-i)$ for $i=0,\ldots,k$.
Since $l_{i}=(t_{0}i+qk_{0})r$,
$l_{i}$ can be divided by $r$ for any $i$.
Especially, we can verify that 
$\gcd\{l_{0},l_{1}\}=r$.
Therefore $\gcd\{l_{0},l_{1},\ldots,l_{k}\}=r$.
In consequence,
there exist $m_{1},\ldots,m_{2^{k}}\in {\bf N}$
such that $\gcd\{m_{1},\ldots,m_{2^{k}}\}=1$ and 
$\ba^{\boxtimes k}=(y^{m_{1}},\ldots,y^{m_{2^{k}}})$ for $y\equiv x^{r}$.
This implies the statement.

As is the case with $p-q>0$,
the case $p-q<0$ can be proved.
\qedh

%
%
\begin{ex}
\label{ex:power}
{\rm
\begin{enumerate}
\item
If $p\in {\bf N}$ and $x\in {\bf R}_{+}$ satisfy
\[x^{p+1}+x^{p}-1=0,\]
then $\varpi_{(x^{p+1},x^{p})}^{\ptimes k}$ is of type ${\rm III}_{x}$
for each $k\in {\bf N}$.
Especially,
when $p=1$,
$\{\varpi_{((\frac{\sqrt{5}-1}{2})^{2},\frac{\sqrt{5}-1}{2})}\}^{\ptimes k}$ 
is of type ${\rm III}_{\frac{\sqrt{5}-1}{2}}$ for each $k\geq 1$.
\item
Let $x$ be a (unique) positive solution of the equation $x^{3}+x-1=0$.
If $\ba=(x,x^{3})$, then 
\[\varpi_{\ba}^{\ptimes k}
\mbox{ is }\left\{
\begin{array}{ll}
\mbox{ {\rm of type} ${\rm III}_{x}$ }\quad & (\mbox{{\rm $k$ is odd}}),\\
\\
\mbox{ {\rm of type} ${\rm III}_{x^{2}}$ }\quad & (\mbox{{\rm $k$ is even}}).\\
\end{array}
\right.
\]
\item
Let $x$ be the positive real solution
of the equation
\[x^{11}+x^{5}-1=0.\]
Define
$\ba=(x^{11},x^{5})$.
Then
\[\varpi_{\ba}^{\ptimes k}
\mbox{ is }\left\{
\begin{array}{ll}
\mbox{ {\rm of type} ${\rm III}_{x^{6}}$ }\quad & 
(\mbox{{\rm $k\equiv 0$ mod $6$}}),\\
\\
\mbox{ {\rm of type} ${\rm III}_{x^{3}}$ }\quad & 
(\mbox{{\rm $k\equiv 3$ mod $6$}}),\\
\\
\mbox{ {\rm of type} ${\rm III}_{x^{2}}$ }\quad & 
(\mbox{{\rm $k\equiv 2$ mod $6$}}),\\
\\
\mbox{ {\rm of type} ${\rm III}_{x}$ }\quad & (\mbox{{\rm otherwise}}).\\
\end{array}
\right.
\]
\end{enumerate}
}
\end{ex}

In consequence,
the tensor power has periodicity arising from
the periodicity of the function ${\bf N}\ni k\mapsto \gcd\{a,k\}$
when $a\geq 1$.\\

\ssfr{Acknowledgments}
The author would like to express his sincere thanks to 
Izumi Ojima for suggestions of the study of 
type III representations of Cuntz-Krieger algebras on many occasions.
He also thanks to Akira Asada 
for a comment of the type of a factor.

\appendix
\section*{Appendix}
%
%
\sftt{Tensor products of type {\rm III} factors}
\label{section:appone}
We review tensor products of type III factors.
In general 
the tensor product of two type III factors
is also a type III factor.
Therefore the tensor product of 
a type ${\rm III}_{\lambda}$ factor
and a type ${\rm III}_{\mu}$ factor
for $0\leq \lambda,\mu\leq 1$
is also a type ${\rm III}_{\nu}$ factor
for a certain $0\leq \nu\leq 1$.
We roughly write this as follows:
%
%
\begin{equation}
\label{eqn:tensorfactor}
{\rm III}_{\lambda}\otimes {\rm III}_{\mu}={\rm III}_{\nu}
\end{equation}
However such $\nu$ is not uniquely determined from
given $\lambda,\mu$ in general \cite{Connesbook}.
Consequently,
(\ref{eqn:tensorfactor})
makes no sense and 
the parameter $0\leq \lambda\leq 1$
of Connes' classification of type III factors
is not compatible with the tensor product of factors in general.

Exceptionally, 
when the L.H.S. in (\ref{eqn:tensorfactor})
is the tensor product of AFD type {\rm III} factors,
$\nu$ is uniquely determined from
given $0<\lambda,\mu\leq 1$ in (\ref{eqn:tensorfactor})
as follows (\cite{Blackadar2006}, III.3.1.14):
\[
\nu=\left\{
\begin{array}{ll}\tau
\quad&\left(
\begin{array}{l}
0<\lambda,\mu<1, \mbox{ and }
\mbox{there exists }(p,q)\in {\bf N}^{2}
\mbox{ such that }\\
 \,\gcd\{p,q\}=1 \mbox{ and }
(\lambda,\mu)=(\tau^{p},\tau^{q})
\end{array}
\right),\\
\\
1 \quad &\mbox{{\rm (otherwise)}}. 
\end{array}
\right.
\]

%
%
\sftt{Several kinds of tensor product of representations}
\label{section:apptwo}
We review basic facts of tensor product of representations.
First, we remark that
there are non-negligible differences between group theory and algebra theory
with respect to terminologies
of tensor product of representations.
In order to explain this, we start with 
tensor products of representations in group theory.

For a representation $\pi_{i}$ of a group $G_{i}$
for $i=1,2$,
define the representation 
$\pi_{1}\otimes_{out} \pi_{2}$ of $G_{1}\times G_{2}$ by
$(\pi_{1}\otimes_{out} \pi_{2})(g_{1},g_{2})\equiv 
\pi_{1}(g_{1})\otimes \pi_{2}(g_{2})$
for $(g_{1},g_{2})\in G_{1}\times G_{2}$.
The representation $\pi_{1}\otimes_{out} \pi_{2}$  
is called the {\it outer tensor product of
representations of $\pi_{1}$ and $\pi_{2}$} \cite{Varadarajan}.
In addition,
when $G_{1}=G_{2}=G$,
let $\iota$ denote the diagonal embedding of $G$ into $G\times G$.
Then 
\[\pi_{1}\otimes_{inn} \pi_{2}\equiv (\pi_{1}\otimes_{out} \pi_{2})\circ\iota\]
is called the 
{\it inner tensor product (or Kronecker product \cite{Tatsuuma}) of $\pi_{1}$ and $\pi_{2}$},
that is, $(\pi_{1}\otimes_{inn} \pi_{2})(g)=
\pi_{1}(g)\otimes_{inn} \pi_{2}(g)$ for $g\in G$.
In usual, 
$\pi_{1}\otimes_{inn} \pi_{2}$ is written as
$\pi_{1}\otimes \pi_{2}$ 
and the inner tensor product of representations is called 
the {\it tensor product of representations} 
for simplicity of description.
Almost always,  
the tensor product of representations in group theory 
means the inner tensor product of representations.

On the other hand,
for a representation $\pi_{i}$ 
of an algebras ${\cal A}_{i}$ for $i=1,2$,
the tensor product representation $\pi_{1}\otimes \pi_{2}$
of $\pi_{1}$ and $\pi_{2}$ is uniquely
defined as the outer tensor product representation 
of $\pi_{1}$ and $\pi_{2}$ for ${\cal A}_{1}\otimes {\cal A}_{2}$
because there is no way to define the inner tensor product
for algebras in general.
Especially, 
the tensor product of representations means
the outer tensor product representation 
in the theory of operator algebras in usual.

If an algebra ${\cal A}$ has a (coassociative) comultiplication 
$\Delta\in {\rm Hom}({\cal A},{\cal A}\otimes {\cal A})$ \cite{Kassel},
then we can define  {\it the (associative) inner tensor product 
(or the Kronecker product \cite{ES})
of representations $\pi_{1}$ and $\pi_{2}$} of ${\cal A}$ by
\[\pi_{1}\otimes_{inn} \pi_{2}\equiv 
(\pi_{1}\otimes \pi_{2})\circ \Delta.\]
Remark that the operation $\otimes_{inn}$
depends on the choice of $\Delta$ in this case.
Here we omit topological problems of tensor product 
of topological algebras for simplification.

At the last,
we see that 
the product $\boxtimes$ of vectors in (\ref{eqn:vector})
can be regarded as the ``Kronecker product of diagonal matrices".
From Theorem \ref{Thm:state},
the ``Kronecker product $\boxtimes$ of vectors" and 
the ``Kronecker product $\ptimes$ of representations" are very close.
This is an interesting accidental coincidence of the terminology
``Kronecker product."

%
%

%

\begin{thebibliography}{99}
%
\bibitem{BP}
A. Berman, R.J. Plemmons,  
Nonnegative matrices in mathematical sciences, 
SIAM, 1994.
%
\bibitem{Blackadar2006}B.  Blackadar, 
Operator algebras. Theory of C$^{*}$-algebras and
von Neumann algebras, Springer-Verlag Berlin Heidelberg New York, 2006.
%
\bibitem{BR2}O. Bratteli,  D.W. Robinson,   
 Operator algebras and quantum statistical mechanics 2,
Springer-Verlag New York, 1981.
%
\bibitem{Connes}A. Connes,  
Une classification des facteurs de type III,
Ann. Sci. \'{E}cole Norm. Sup. (4), 6 (1973) 133--252.
%
\bibitem{Connesbook}A.  Connes, 
Non commutative geometry, Academic Press, Orlando, 1993.
%
%
\bibitem{CK}J. Cuntz, W. Krieger,  
A class of $C^{*}$-algebra and topological Markov chains,
Invent. Math. 56 (1980) 251--268.
%
%
\bibitem{Dief}A.S. Dief,  
Advanced matrix theory for scientists and engineers,
Abacus Press, Turnbridge Wells \& London, 1982.
%
\bibitem{Dixmier}J. Dixmier,
$C^{*}$-algebras,
North-Holland Publishing Company
Amsterdam New York Oxford, 1977.
%
\bibitem{ES}M. Enock,   J.M. Schwartz, 
Kac algebras and duality of locally compact groups,
Springer-Verlag, 1992.
%
\bibitem{EFW}M. Enomoto,  M. Fujii, Y. Watatani,  
KMS states for gauge action on $\coa$, 
Math. Japon. 29(4) (1984) 607--619.
%
\bibitem{Evans}D.E. Evans,  
On ${\cal O}_{n}$, 
Publ. Res. Inst. Math. Sci. 16 (1980) 915--927.
%
\bibitem{EL}R. Exel,  M. Laca,  
Cuntz-Krieger algebras for infinite matrices, 
J. reine angew. Math.(Crelle) 512 (1999) 119--172.
%
\bibitem{EL2}R. Exel,  M. Laca,  
Partial dynamical systems and the KMS condition, 
Commun. Math. Phys. 232(2) (2003) 223--277. 
%
\bibitem{Izumi}M. Izumi,   
Subalgebras of infinite C$^{*}$-algebras with 
finite Watatani indices. I. Cuntz algebras,
Commun. Math. Phys. 155(1) (1993) 157--182.
%
%
\bibitem{Kassel}C. Kassel, 
Quantum groups, Springer-Verlag, 1995.
%
\bibitem{CK01}K. Kawamura,  
Polynomial embedding of Cuntz-Krieger algebra into Cuntz algebra,
preprint  RIMS-1391. 
%
\bibitem{CKR01}K. Kawamura, 
Representations of the Cuntz-Krieger algebras. I ---General theory---,
preprint RIMS-1454.
%
\bibitem{CKR02}K. Kawamura, 
Representations of the Cuntz-Krieger algebras. II
---Permutative representations---,
 preprint RIMS-1462.
%
\bibitem{PFO01}K. Kawamura, 
The Perron-Frobenius operators, invariant measures 
and representations of the Cuntz-Krieger algebras,
J. Math. Phys. 46(8) (2005) 083514-1--083514-6. 
%
\bibitem{PE01}K.  Kawamura,  
Polynomial endomorphisms of  the Cuntz algebras arising from
permutations. I ---General theory---,
Lett. Math. Phys. 71 (2005) 149--158. 
%
\bibitem{PE02}K. Kawamura,   
Branching laws for polynomial endomorphisms of Cuntz algebras
arising from permutations,
Lett. Math. Phys. 77 (2006) 111--126.
%
\bibitem{PE03}K. Kawamura,  
Automata computation of branching laws for endomorphisms of 
Cuntz algebras,
Int. J. Alg. Comput. 17(7) (2007) 1389--1409.
%
\bibitem{TS01}K.  Kawamura, 
A tensor product of representations of Cuntz algebras,
Lett. Math. Phys. 82 (2007) 91--104.
%
\bibitem{TS02}K. Kawamura,  
C$^{*}$-bialgebra defined by the direct sum of Cuntz algebras,
J. Algebra,  319 (2008) 3935--3959.
%
\bibitem{TS05}K. Kawamura, 
C$^{*}$-bialgebra defined by the direct sum of Cuntz-Krieger algebras,
math.OA/0801.4597v2.
%
\bibitem{Okayasu}R. Okayasu,  
Type {\rm III} factors arising from Cuntz-Krieger algebras,
Proc. Amer. Math. Soc. 131(7) (2002) 2145--2153.
%
\bibitem{Ped}G.K. Pedersen, 
$C^{*}$-algebras and their automorphism groups,
Academic Press, 1979.
%
\bibitem{RS}M. R\o rdam,  E. St\o rmer, 
Classification of nuclear C$^{*}$-algebras.
Entropy in operator algebras,
Springer-Verlag Berlin Heidelberg, 2002.
%
\bibitem{Seneta}E. Seneta,  
Non-negative matrices and Markov chains. 2nd. edit.,
Springer-Verlag New York Heidelberg Berlin, 1981.
%
%
%
\bibitem{Tatsuuma}N. Tatsuuma,
Decomposition of Kronecker products of 
representations of the inhomogeneous Lorentz group,
Proc. Japan Acad. 38 (1962) 156--160.
%
\bibitem{Varadarajan}V.S. Varadarajan,
Lie groups, Lie algebras, and their representations,
Springer-Verlag New York Berlin Heidelberg Tokyo, 1984.
%
\bibitem{Wallach}N.R. Wallach,  
Real reductive groups II,
Academic Press, 1992.
%
\end{thebibliography}
\end{document}